%&amstex          
\input amstex\documentstyle{amsppt}  
\pagewidth{12.5cm}\pageheight{19cm}\magnification\magstep1
\topmatter
\title{Families of isotropic subspaces in a symplectic $\bold Z/2$-vector space}
\endtitle
\author G. Lusztig\endauthor
\address{Department of Mathematics, M.I.T., Cambridge, MA 02139}\endaddress
\thanks{Supported by NSF grant DMS-2153741}\endthanks
\endtopmatter   
\document

\define\bV{\bar V}

\define\mpb{\medpagebreak}

\define\be{\bar e}

\define\hs{\hat s}

\define\pe{\perp}
\define\si{\sim}

\define\sqc{\sqcup}

\define\lb{\linebreak}

\define\op{\oplus}
   
\define\part{\partial}
\define\emp{\emptyset}

\define\n{\notin}

\define\m{\mapsto}
\define\do{\dots}

\define\sub{\subset}    

\define\T{\times}
\define\ti{\tilde}
\define\nl{\newline}
\redefine\i{^{-1}}

\define\un{\underline}
\define\ov{\overline}

\define\supp{\text{\rm supp}}

\redefine\spa{\spadesuit}

\redefine\b{\beta}
\redefine\c{\chi}

\redefine\d{\delta}
\define\e{\epsilon}

\define\io{\iota}
\redefine\o{\omega}
\define\p{\pi}
\define\ph{\phi}

\define\s{\sigma}
\redefine\t{\tau}

\define\z{\zeta}
\define\x{\xi}

\define\ee{\bold e}

\define\CC{\bold C}

\define\NN{\bold N}

\define\ZZ{\bold Z}

\define\cf{\Cal F}

\define\ci{\Cal I}

\define\cx{\Cal X}

\define\fE{\frak E}

\define\ty{\ti y}

\define\tB{\ti B}

\define\bp{\bar p}

\head Introduction\endhead
\subhead 0.1 \endsubhead
Let $F=\ZZ/2$ be the field with two elements.
Let $\bV$ be an $F$-vector space of finite dimension $2n\ge2$
endowed with a nondegenerate symplectic form \lb $<,>$ and with
a collection of vectors $\be_0,\be_1,\be_2,\do,\be_{2n}$ such that 

$<\be_0,\be_1>=<\be_1,\be_2>=\do=<\be_{2n-1},\be_{2n}>=
<\be_{2n},\be_0>=1$,

$<\be_1,\be_0>=<\be_2,\be_1>=\do=<\be_{2n},\be_{2n-1}>=
<\be_0,\be_{2n}=1$

and $<\be_i,\be_j>=0$ for all other pairs $i,j$. 
(Such a collection is called a ``circular basis'' in \cite{L20a}.)

In \cite{L20a} we have introduced a family $\cf(\bV)$ of isotropic
subspaces of $\bV$ with remarkable properties:

{\it There is a unique bijection $\cf(\bV)@>\si>>\bV$ such that
any $x\in\bV$ is contained in the corresponding subspace of $\bV$. The
characteristic functions of the various subspaces in $\cf(\bV)$ form a
new basis of the complex vector space $\bV^\CC$ of functions
$\bV@>>>\CC$ which is related to the obvious basis of $\bV^\CC$ by an
upper triangular matrix with $1$ on diagonal (in some partial order
$\le$ on $\cf(\bV)$).}

(In fact the collection $\cf(\bV)$ was already introduced 
in \cite{L20}, but in a less symmetric form.)

A further property of $\cf(\bV)$ was found in \cite{L20a}, namely that
the matrix of the Fourier transform $\bV^\CC@>>>\bV^\CC$ with respect to
the new basis is upper triangular with $\pm1$ on diagonal. The proof of
this property was based on the observation that the new basis admits a
dihedral symmetry which was not visible in the definition of
\cite{L20}. 

In this paper we give a new non-inductive definition of $\cf(\bV)$
which is visibly compatible with
the dihedral symmetry (the definition of \cite{L20} has no such a
symmetry property; the definition in \cite{L20a} did have the symmetry
property but was inductive). We also give a formula for the bijection
$\cf(\bV)@>\si>>\bV$ above which is clearly compatible
with the dihedral symmetry. (See Theorem 1.4.)

Let $V$ be an $F$-vector space with basis $e_0,e_1,\do,e_{2n}$
such that $\bV$ is the quotient of $V$ by the line
$F(e_0+e_1+\do+e_{2n})$ and $\be_i$ is the image of $e_i$ under
the obvious map $V@>>>\bV$. In \S4 we define an analogue $\ti\cf(V)$
of $\cf(\bV)$ which is a refinement of $\cf(\bV)$ and has several
properties of $\cf(\bV)$.

In \S5-\S7 we study a modification of the family $\cf(\bV)$ which plays
the same role in the theory of unipotent representations of orthogonal
groups over a finite field as that played by $\cf(\bV)$ in the
analogous theory for symplectic groups over a finite field.

\head 1. Statement of the Theorem\endhead
\subhead 1.1\endsubhead
Let $V$ be an $F$-vector space endowed with a symplectic form
$<,>:V\T V@>>>F$ and a map $e:S@>>>V$, $s\m e_s$ where S is a finite
set. Let $\fE$ be the set of unordered pairs $s\ne s'$ in $S$ such that
$<e_s,e_{s'}>=1$. This is the set of edges of a graph with set of
vertices $S$.
For any $I\sub S$ we set $e_I=\sum_{s\in I}e_s\in V$ and we denote by
$\un I$ the full subgraph of $(S,\fE)$ whose set of vertices is $I$.
Let $\ci$ be the set of all $I\sub S$ such that $\un I$
is a graph of type $A_m$ for some $m\ge1$. We have $\ci=\ci^0\sqc\ci^1$
where $\ci^0=\{I\in\ci;|I|=0\mod2\}$, $\ci^1=\{I\in\ci;|I|=1\mod2\}$.
For $I,I'$ in $\ci^1$ we write $I\prec I'$ whenever
$I\subsetneqq I'$ and
$\un{I'-I}$ is disconnected.
For $I,I'$ in $\ci^1$ we write $I\spa I'$ whenever $I\cap I'=\emp$ and
$\un{I\cup I'}$ is disconnected.
For $I\in\ci^1$ let $I^{ev}$ be the set of all $s\in I$ such that
$I-\{s\}=I'\sqc I''$, with $I'\in\ci^1$, $I''\in\ci^1$, $I'\spa I''$.
Let $I^{odd}=I-I^{ev}$. We have $|I^{ev}|=(|I|-1)/2$. 

\subhead 1.2\endsubhead
Let $R$ be the set whose elements are finite unordered sequences of
objects of $\ci^1$. For $B\in R$ let $L_B$ be the subspace of $V$
generated by $\{e_I;I\in B\}$; for a subspace $L$ of $V$ let
$B_L=\{I\in\ci^1;e_I\in L\}\sub R$. For $s\in S$, $B\in R$ we set
$$g_s(B)=|\{I\in B;s\in I\}|$$
(here $|?|$ denotes the number of elements of $?$) and
$$\e_s(B)=(1/2)g_s(B)(g_s(B)+1)\in F.$$
For $B\in R$ we set
$$\e(B)=\sum_{s\in S}\e_s(B)e_s\in V.$$
For $B\in R$ we set $\supp(B)=\cup_{I\in B}I\sub S$.

Let $\ph(V)$ be the set consisting of all $B\in R$ such that
$(P_0),(P_1)$ below hold.

$(P_0)$ If $I\in B,I'\in B$, then $I=I'$, or $I\spa I'$, or
$I\prec I'$, or $I'\prec I$.

$(P_1)$ Let $I\in B$. There exist $I_1,I_2,\do,I_k$ in $B$ such that
$I^{ev}\sub I_1\cup I_2\cup\do\cup I_k$ (disjoint union),
$I_1\prec I$, $I_2\prec I,\do,I_k\prec I$.
\nl
We say that $(V,<,>,e)$ is {\it perfect} if properties (i)-(iv) below
hold.

(i) If $B\in\ph(V)$, then $\{e_I;I\in B\}$ is a basis of $L:=L_B$;
moreover, $B=B_L$.

(ii) For any $B\in\ph(V)$ we have $\e(B)\in L_B$. Hence $\e$ restricts
to a map $\ph(V)@>>>V_0$ (denoted again by $\e$) where
$V_0=\cup_{B\in\ph(V)}L_B\sub V$.

(iii) The map $\e:\ph(V)@>>>V_0$ is a bijection. 

(iv) If $B,B'$ in $\ph(V)$ are such that $\e(B')\in L_B$, then
$g_s(B')\le g_s(B)$ for any $s\in S$.

\mpb

For $B',B$ in $\ph(V)$ we say that $B'\le B$ if there exist
$B_0,B_1,B_2,\do,B_k$ in $\ph(V)$ such that $B_0=B',B_k=B$,

$\e(B_0)\in L_{B_1},\e(B_1)\in L_{B_2},\do,\e(B_{k-1})\in L_{B_k}$.
\nl
We show:

(a) {\it If $(V,<,>,e)$ is  perfect, then $\le$ is a partial order on
$\ph(V)$.}
\nl
Assume that we have elements $B_0,B_1,\do,B_k,B'_0,B'_1,\do,B'_l$ in
$\ph(V)$ such that 

$\e(B_0)\in L_{B_1},\e(B_1)\in L_{B_2},\do,\e(B_{k-1})\in L_{B_k}$,

$\e(B'_0)\in L_{B'_1},\e(B'_1)\in L_{B'_2},\do,
\e(B'_{l-1})\in L_{B'_l}$,
\nl
and $B_0=B'_l,B'_0=B_k$. We must prove that $B_0=B'_0$
Using (iv) and our assumptions we have for any $s\in S$:

$g_s(B_0)\le g_s(B_1)\le g_s(B_2)\le\do\le g_s(B_k)=g_s(B'_0)$,

$g_s(B'_0)\le g_s(B'_1)\le g_s(B'_2)\le\do\le g_s(B'_l)=g_s(B_0)$.
\nl
It follows that $g_s(B_0)\le g_s(B'_0), g_s(B'_0)\le g_s(B_0)$, so that
$g_s(B_0)=g_s(B'_0)$. Since this holds for any $s$, we see that
$\e(B_0)=\e(B'_0)$. Using the injectivity of $\e$ (see (iii)), we
deduce that $B_0=B'_0$, as desired.

\subhead 1.3 \endsubhead
We will consider three cases:

(a) $V,<,>,e:S@>>>V$ are such that $\{e_s;s\in S\}$ is a basis of $V$ and
$(S,\fE)$ is a graph of type $A_{N-1}$, $N\in\{3,5,7,\do\}$;

(b) $V,<,>,e:S@>>>V$ are such that $\{e_s;s\in S\}$ is a basis of $V$ and
$(S,\fE)$ is a graph of affine type $A_{N-1}$, $N\in\{3,5,7,\do\}$;

(c) $V,<,>,e:S@>>>V$ in (b) are replaced by $\bV=V/Fe_S$, by the
symplectic form induced by $<,>$ (denoted again by $<,>$), and by
$\p e:S@>>>\bV$, where $\p:V@>>>\bV$ is the obvious map.

\mpb

In cases (b),(c) we note that the automorphism group of the graph
$(S,\fE)$ is a dihedral group $Di_{2N}$ of order $2N$. It acts
naturally on $V$ in (b) by permutations of the basis; this induces
an action of $Di_{2N}$ on $\bV$ in (c).

Let $I\sub S$; in cases (b),(c) we assume that $I\ne S$.
There is a well defined subset $c(I)$ of
$\ci$ such that $I'\spa I''$ for any $I'\ne I''$ in $c(I)$ and
$I=\sqc_{I'\in c(I)}I'$. Note that $\{\un I';I'\in c(I)\}$ are the
connected components of the graph $\un I$.

We now state the following result.
\proclaim{Theorem 1.4}In each of the cases 1.3(a),(b),(c), $(V,<,>,e)$
is perfect.
\endproclaim

\subhead 1.5\endsubhead
In case 1.3(a), Theorem 1.4 is contained in \cite{L19}. Let $\cf(V)$
be the set of subspaces of $V$ of the form $L_B$ for some $B\in\ph(V)$.
Note that $B\m L_B$ is a bijection $\ph(V)@>\si>>\cf(V)$.

We can write the elements of $S$ as a sequence $s_1,s_2,\do,s_{N-1}$
in which any two consecutive elements are joined in the graph $(S,\fE)$.
Let $I\sub S$. Let $c(I)$ be as in 1.3.
Let $c(I)^{0+}$ (resp. $c(I)^{0-}$) be the set of all $I'\in c(I)$
such that $I'=\{s_k,s_{k+1},\do,s_l\}$ where $k$ is even, $l$ is odd
(resp. $k$ is odd, $l$ is even). Let $V_0$ be the subset of $V$
consisting of all $e_I$ where $I\sub S$ satisfies
$|c(I)^{0+}|=|c(I)^{0-}|$. From \cite{L19} it is known that $V_0$
coincides with the subset of $V$ appearing in 1.2(ii) that is,
n
(a) $\cup_{L\in\cf(V)}L=V_0$.

\head 2. The case 1.3(c)\endhead
\subhead 2.1\endsubhead
In this section we assume that we are in case 1.3(c).
For $s\in S$ we set $\be_s=\p(e(s))$. For $I\sub S$ we set
$\be_I=\sum_{s\in I}\be_s$.
Note that $\{\be_s;s\in S\}$ is a circular basis of $\bV$
(in the sense of \cite{L20a}) and to this we can attach a collection
$\cf(\bV)$ of subspaces of $\bV$ as in \cite{L20a}. We recall how this
was done. For any $s\in S$ we set
$$\hs=\{s'\in S;<\be_s,\be_{s'}>=1\}\cup\{s\}\sub S.$$
We have $|\hs|=3$. We set $\be_s^\pe=\{x\in\bV;<x,\be_s>=0\}$ and
$\bV_s=\be_s^\pe/F\be_s$. This is a symplectic $F$-vector space with
circular basis $\{\be_{s'};s'\in S-\hs\}\sqc\{\be_{\hs}\}$.
Thus the analogue of $S$ when $\bV$ is replaced by $\bV_s$ is
$S_s=(S-\hs)\sqc\{\hs\}$ (a set with $|S|-2$ elements).
Let $\bp_s:\be_s^\pe@>>>\bV_s$ be the obvious linear map.
We define a collection $\cf(\bV)$ of subspaces of $\bV$ by
induction on $N$. If $N=3$, $\cf(\bV)$ consists of $0$ and of
$\bp_s\i(0)$ for various $s\in S$. If $N\ge5$, $\cf(\bV)$ consists of
$0$ and of $\bp_s\i(L')$ for various $s\in S$ and various
$L'\in\cf(\bV_s)$ (which is defined by the induction hypothesis).
In \cite{L20a}, $\cf(\bV)$ is also identified with a collection
of subspaces of $\bV$ introduced in \cite{L20} in terms of a chosen
element $t\in S$. From this identification we see that:

(a) {\it if $L\in\cf(\bV)$ and
$B_L^t:=\{I\in\ci;I\sub S-\{t\},\be_I\in L\}$,
then $\{\be_I;I\in B_L^t\}$ is an $F$-basis of $L$, so that
$L=L_{B_L^t}$.}
\nl
Now if $I\in\ci$, then $S-I\in\ci$ and we have $\be_I=\be_{S-I}$.
Moreover, exactly one of $I,S-I$ is contained in $S-\{t\}$ and exactly
one of $I,S-I$ is in $\ci^1$. We deduce that:

(b) {\it If $L\in\cf(\bV)$, and
$$B_L:=\{I\in\ci^1;\be_I\in L\}=
\{I\in\ci^1;I\in B_L^t\}\sqc\{I\in\ci^1;S-I\in B_L^t\}$$
then $\{\be_I;I\in B_L\}$ is an $F$-basis of $L$, so that $L=L_{B_L}$.}

\subhead 2.2\endsubhead
We show that for $B\in R$:

(a) {\it we have $B\in\ph(\bV)$ if and only if $L_B\in\cf(\bV)$.}
\nl
The proof is analogous to that of the similar result in case 1.3(a)
given in \cite{L19}. We argue by induction on $N$. If $N=3$, (a) is
easily verified. In this case, $B$ is either $\emp$ or it is of the
form $\{s\}$ for some $s\in S$. We now assume that $N\ge5$. For
$s\in S$ we denote by $\ci^1_s,R_s$ the analogues of $\ci^1,R$ when
$S$ is replaced by $S_s$ (see 2.1). For $J\in\ci^1_s$ we write
$\be_J\in\bV_s$ for the analogue of
$\be_I\in\bV, I\in\ci^1$. We have
$$\bp_s\i(\be_J)=\{\be_I,\be_I+\be_s\}$$
for a well defined $I\in\ci^1$ such that $s\n I$; we set
$I=\x_s(J)$. There is a well defined map $\t_s:R_s@>>>R$, $B'_1\m B_1$
where $B_1$ consists of $\{s\}$ and of all $\x_s(J)$ with $J\in B'_1$.
From the definitions we see that (assuming that $B'_1\in R_s$ and
$B_1=\t_s(B'_1)$), the following holds.

(b) {\it $B'_1$ satisfies $(P_0)$ if and only if $B_1$ satisfies
$(P_0)$; $B'_1$ satisfies $(P_1)$ if and only if $B_1$ satisfies
$(P_1)$.}
\nl
Assume now that $B$ is such that $L:=L_B\in\cf(\bV)$, so
that $B=B_L$. We
show that $B$ satisfies $(P_0),(P_1)$. If $B=\emp$, this is obvious.
If $B\ne\emp$, we have $L=\bp_s\i(L')$ where $s\in S,L'\in\cf(\bV_s)$.
From the definition we have $\t_s(B_{L'})=B_L$. By the induction
hypothesis, $B_{L'}$ satisfies $(P_0),(P_1)$; using (b), we see that
$B=B_L$ satisfies $(P_0),(P_1)$.

Conversely, assume that $B$ satisfies $(P_0),(P_1)$. We show that
$B=B_L$ for some $L\in\cf(\bV)$. If $B=\emp$ this is obvious. Thus we
can assume that $B\ne\emp$. Let $I\in B$ be such that $|I|$ is minimum.
If $s\in I^{ev}$ (see 1.1) then by $(P_1)$ we can find $I'\in B$ with
$s\in I',|I'|<|I|$, a contradiction. We see that $I^{ev}=\emp$. Thus,
$I=\{s\}$ for some $s\in S$. Using $(P_0)$ and $\{s\}\in B$, we see
that for any $I'\in B-\{s\}$ we have $\{s\}\prec I'$ or $I'\spa\{s\}$.
It follows that $B=\t_s(B')$ for some $B'\in R_s$. From (b) we see that
$B'$ satisfies $(P_0),(P_1)$. From the induction hypothesis we see that
$B'=B_{L'}$ for some $L'\in\cf(\bV_s)$. Let $L=\bp_s\i(L')$. We have
$L\in\cf(\bV)$ and $B=B_L$. This proves (a).

\mpb

We see that we have a bijection

(c) $\ph(\bV)@>\si>>\cf(\bV)$, $B\m L_B$.
\nl
Using now 2.1(b) we see that 1.2(i) holds for any $B\in\ph(\bV)$.

\subhead 2.3\endsubhead
We now fix $t\in S$.
Let $B\in\cf(\bV)$, let $L=L_B\in\cf(\bV)$ and let $B^t=B^t_L$ (see 2.1).
For any $s\in S-\{t\}$ we set
$$f_s(B)=|\{I\in B^t\cap\ci^1;s\in I\}|
-|\{I\in B^t\cap\ci^0;s\in I\}|-\un{|B^t\cap\ci^0|}$$
where for any $m\in\ZZ$ we set $\un m=0$
if $m$ is even, $\un m=1$ if $m$ is odd. We also set
$$\e'(B)=\sum_{s\in S-\{t\}}(1/2)f_s(B)(f_s(B)+1)\be_s\in\bV.$$
From \cite{L20},\cite{L20a} we see using 2.2(c) that:

(a) {\it we have $\e'(B)\in L_B$ for any $B\in\ph(\bV)$ and
$B\m\e'(B)$ defines a bijection
$\e':\ph(\bV)@>\si>>\bV$.}

\subhead 2.4\endsubhead
We wish to rewrite the bijection $\e':\ph(\bV)@>\si>>\bV$ without
reference to $t\in S$. Recall that for any
$B\in\ph(\bV)$ and any $s\in S$ we have

(a) $g_s(B)=|\{I\in B;s\in I\}|\in\NN$.
\nl
Setting $\b=|B^t\cap\ci^0|$ where
$B^t=B^t_L$, $L=L_B$ (see 2.1) we have

(b) $g_t(B)=\b$.
\nl
For $s\in S-\{t\}$ we show:

(c) $f_s(B)=g_s(B)-\b-\un\b$ 
\nl
that is,
$$|\{I\in B^t\cap\ci^1;s\in I\}|-|\{I\in B^t\cap\ci^0;s\in I\}|
=|\{I\in B;s\in I\}|-\b.$$
To prove this, we substitute $|\{I\in B;s\in I\}|$ by
$$|\{I\in B^t\cap\ci^1;s\in I\}|+|\{I\in B^t\cap\ci^0;s\n I\}|.$$
We see that desired equality becomes
$$\align&|\{I\in B^t\cap\ci^1;s\in I\}|-|\{I\in B^t\cap\ci^0;s\in I\}|
\\&=|\{I\in B^t\cap\ci^1;s\in I\}|+|\{I\in B^t\cap\ci^0;s\n I\}|-\b
\endalign$$
which is obvious.

We shall prove the following formula for $\e'(B)$:
$$\e'(B)=\sum_{s\in S}(1/2)g_s(B)(g_s(B)+1)\be_s\tag d.$$
Using (c) we have for $s\in S-\{t\}$:
$$\align&(1/2)f_s(B)(f_s(B)+1)=(1/2)(g_s(B)-\b-\un\b)(g_s(B)-\b-\un\b+1)
\\&=(1/2)g_s(B)(g_s(B)+1)+H\endalign$$
where
$$H=(1/2)(g_s(B)(-2\b-2\un\b)+(\b+\un\b)^2-\b-\un\b).$$
Note that
$$-2\b-2\un\b=0\mod4, (\b+\un\b)^2=0\mod4, -\b-\un\b=-\b(\b+1)\mod4$$
hence $H=-\b(\b+1)\mod2$. Thus,
$$\align&\e'(B)=\sum_{s\in S-\{t\}}(1/2)g_s(B)(g_s(B)+1)\be_s+\\&
\sum_{s\in S-\{t\}}(1/2)g_t(B)(g_t(B)+1)\be_s=
\sum_{s\in S}(1/2)g_s(B)(g_s(B)+1)\be_s.\endalign$$
We have used that $\sum_{s\in S}\be_s=0$. This proves (d).

From (d) and 2.3(a) we see that 1.2(ii),(iii) hold in our case with
$\bV_0=\bV$; moreover, $\e'$ in 2.3 is the same as $\e$ in 1.2.

\subhead 2.5\endsubhead
From the results in \cite{L20},\cite{L20a} it is known that if
$B,B'$ in $\ph(\bV)$ satisfy $\e'(B')\in L_B$
(that is, $\e(B')\in L_B$), then $f_s(B')\le f_s(B)$
for any $s\in S-\{t\}$ and $|B^t_{L'}\cap\ci^0|\le|B^t_L\cap\ci^0|$.
(Notation of 2.1 with $L=L_B,L'=L_{B'}$.) We show that 

(a) {\it $g_s(B')\le g_s(B)$ for any $s\in S$.}
 \nl
When $s=t$ this follows from  2.4(b). We now assume that $s\ne t$.
Using 2.4(c) we have
$$g_s(B')+g_t(B')+\un{g_t(B')}\le g_s(B)+g_t(B)+\un{g_t(B)}$$
hence it is enough to show that

(b)  $g_t(B)-g_t(B')+\un{g_t(B)}-\un{g_t(B')}\ge0$.
\nl
If $g_t(B')=g_t(B)$, then (b) is obvious. Assume now that
$g_t(B')\ne g_t(B)$. As we have seen above, we have
$g_t(B')\le g_t(B)$ hence $g_t(B)-g_t(B')\ge1$. We have
$\un{g_t(B)}-\un{g_t(B')}\in\{0,1,-1\}$, hence (b) holds.
This proves (a).

We see that 1.2(iv) holds in our case. Thus Theorem 1.4 is proved
in case 1.3(c).

In the remainder of this paper we write $\bar\e$ instead of
$\e:\ph(\bV)@>>>\bV$ to distinguish it from $\e$ in cases 1.3(a),(b).

\subhead 2.6\endsubhead
We note:

(a) {\it If $B\in\ph(\bV)$, then $\supp(B)\ne S$.}
\nl
This holds since $B$ has property $(P_0)$.

\subhead 2.7\endsubhead
For $t\in S$ let $V(t)$ be the $F$-subspace of $V$ with basis
$\{e_s;s\in S-\{t\}\}$. Then $V(t)$ with this basis and the
restriction of $<,>$ is as in 1.3(a).
Let $R(t)$ be the analogue of $R$ when $V$ in 1.3(a) is replaced by
$V(t)$; we have $R(t)\sub R$. 
Then $\ph(V(t))$ (a collection of elements of $R(t)$) is
defined. From the definition we have $\ph(V(t))\sub\ph(\bV)$.
 Now let $B\in\ph(\bV)$. By 2.6(a) we can find
$t\in S$ such that $\supp(B)\sub S-\{t\}$. Now $B$ satisfies
$(P_0),(P_1)$ relative to $V(t)$. Hence we have $B\in\ph(V(t))$.
We see that

(a) $\ph(\bV)=\cup_{t\in S}\ph(V(t))$.
\nl
From the definitions we see that for any $t\in S$ the following diagram
is commutative:
$$\CD
\ph(V(t))@>>>\ph(\bV)\\
@V\e VV         @V\bar\e VV  \\
V(t)_0@>>>\bV
\endCD$$
Here the left vertical maps are as in 1.2; the horizontal maps are
the obvious inclusions.

\subhead 2.8\endsubhead
We wish to compare the approach to $\ph(\bV)$ given in this paper with that in \cite{L23}. Let $S'=\fE$. We can regard $S'$ as a set of
vertices of a graph in which $\{s_1,s_2\}\in\fE,\{s_3,s_4\}\in\fE$
are joined whenever $|\{s_1,s_2\}\cap\{s_3,s_4\}|=1$. Thus the set
$\fE'$ of edges of this graph is in obvious bijection with $S$.
Note that the graph $(S',\fE')$ is isomorphic to $(S,\fE)$ hence the
analogues $\bV',\ci'{}^1,\ph(\bV')$ of $\bV,\ci^1,\ph(\bV)$ when
$(S,\fE)$ is replaced by $(S',\fE')$ are defined. We can view $\bV'$
as the $F$-vector space consisting of all subsets of $S$ of even
cardinal in which the sum of $X,X'$ is $(X\cup X'(-(X\cap X')$, which
is endowed with the symplectic form $X,X'\m|X\cap X'|\mod2$ and with a
circular basis consisting of all two elements subsets of $S$ which are
in $\fE$. This circular basis is therefore indexed by $S'$.
Now an object of $\ci'{}^1$ is a subgraph of type $A_{2k+1}$ ($k\ge0$)
of $S'$, that is with vertices of the form
$\{s_1,s_2\},\{s_2,s_3\},\do,\{s_{2k+1},s_{2k+2}\}$; this is the same
as a graph of type $A_{2k+2}$ of $S$ (with vertices $s_1,s_2,\do,s_{2k+2}$) and is completely determined by
the pair of (distinct) elements $s_1,s_{2k+2}$. Thus $\ci'{}^1$ can be
identified with the set of two element subsets of $S$. In this way
$\ci'{}^1$ appears as a subset of $\bV'$ and each $X$ in
$\ci'{}^1$ determines a subgraph of type $A_{2k+2}$ ($k\ge0$) of
$S$; the set of vertices of this subgraph is denoted by $\un X$.
(We have $\un X\sub\bV'$ and $X\sub\un X$.)

Now $\ph(\bV')$ becomes the set of all unordered pairs
$X_1,X_2,\do,X_k$ of two element subsets of $S$ such that
$X_i\cap X_j=\emp$ for $i\ne j$ and such that for any
$i\in\{1,2,\do,k\}$ there exists $j_1<j_2<\do<j_s$ in $\{1,2,\do,k\}$
such that
$$\un{X_i}-X_i=
\un{X_{j_1}}\sqc\un{X_{j_2}}\sqc\do\sqc\un{X_{j_s}}.$$

This approach appears in \cite{L23} (in a less symmetric and more
complicated way) where $S$ is taken to be
$S_N=\{1,2,\do,N\}$ with $\fE$ consisting of
$\{1,2\},\{2,3\},\do,\{N-1,N\},\{N,1\}$.

The set $\cx_{N-1}$ defined in \cite{L23, 1.3} is the same as
$\ph(\bV')$ although its definition is less symmetric and
more complicated. Hence it is the same as $\ph(\bV)$ if
$\bV,\bV'$ are identified by $\be_s\m\{s,s+1\}$ if
$s\in\{1,2,\do,N-1\}$ and $\be_N\m\{N,1\}$.

\head 3. The case 1.3(b)\endhead
\subhead 3.1 \endsubhead
In this section we assume that we are in the setup of 1.3(b).
Let $V_0$ be the set of all vectors of $V$ which are of the form 
$e_I$ with $I\sub S,I\ne\emp,I\ne S$ such that $|c(I)\cap\ci^0|$ is
even (here $c(I)\sub\ci$ is as in 1.4); let $V_1$ be the set of all
vectors of $V$ which are of the form $e_S$ or $e_I$ with
$I\sub S,I\ne \emp,I\ne S$ such that $|c(I)\cap\ci^0|$ is odd.
We have clearly:

(a) $V=V_0\sqc V_1$.
\nl
We show:

(b) {\it If $I\sub S, I\ne\emp,I\ne S$, then $e_I\in V_0$ if and
only if $e_{S-I}\in V_1$. In particular, $x\m x+e_S$ is a bijection
$V_0@>\si>>V_1$.}
\nl
We have $c(I)=\{I_1,I_3,\do,I_{2r-1}\}$,
$c(S-I)=\{I_2,I_4,\do,I_{2r}\}$ and (if $r>1$) we have
$I_1\cup I_2\in\ci$, $I_2\cup I_3\in\ci$, $\do$,
$I_{2r-1}\cup I_{2r}\in\ci$, $I_{2r}\cup I_1\in\ci$; in particular,
we have $|c(I)|=|c(S-I)|$. (This remains true also when $r=1$.) Hence,
setting $c^0(I)=c(I)\cap\ci^0$, $c^1(I)=c(I)\cap\ci^1$, we have
$$|c^0(I)|-|c^0(S-I)|=-|c^1(I)|+|c^1(S-I)|.$$
Modulo 2 this equals
$$\align&
|c^1(I)|+|c^1(S-X)|=\sum_{I'\in c^1(I)}|I'|+\sum_{I'\in c^1(S-I)}|I'|
\\&=\sum_{I'\in c^1(I)}|I'|+\sum_{I'\in c^1(S-I)}|I'|+
\sum_{I'\in c^0(I)}|I'|+\sum_{I'\in c^0(S-I)}|I'|\\&
=\sum_{I'\in c(I)}|I'|+\sum_{I\in c(S-I)}|I'|=|I|+|S-I|=|S|.\endalign$$
Since $|S|$ is odd, we see that

(c) $|c^0(I)|-|c^0(S-I)|=1\mod2$
\nl
so that (b) holds.

We show:

(d) {\it Let $\p_0:V_0@>>>\bV$ be the restriction of $\p:V@>>>\bV$.
Then $\p_0$ is a bijection.}
\nl
Assume that $v\ne v'$ in $V_0$ satisfy $\p(v)=\p(v')$. If $v=0$, then
$v'\in\p\i(0)-\{0\}$ hence $v'=e_S$.
But $e_S\n V_0$, a contradiction. If $v\ne0$, then $v=e_I,v'=e_{S-I}$
with $I\sub S,I\ne\emp,I\ne S$. Now $|c^0(I)|$ is even, $|c^0(S-I)|$
is even; but the sum of these numbers is odd by (c), a contradiction.
We see that $\p_0$ is injective.

From (b) we see that $|V_0|=|V_1|$ so that both of these numbers
are equal to $(1/2)|V|=2^{N-1}$. We see that $\p_0$ is an injective map
between two finite sets with $2^{N-1}$ elements;
hence it is a bijection. This proves (d).

\subhead 3.2\endsubhead
Note that the sets $R,\ci$ for this $V$ and for $\bV$ in 1.3(c) are
the same. Hence we have $\ph(V)=\ph(\bV)$.
For $B\in\ph(V)$ we denote by $M_B$ (resp. $L_B$) the subspace of $V$
(resp. $\bV$) generated by $\{e_I;I\in B\}$ (resp. $\{\be_I;I\in B\}$).
Since $\{\be_I;I\in B\}$ is a basis of $L_B$, we see that
$\{e_I;I\in B\}$ is a basis of $M_B$ and that $\p$ restricts to an
isomorphism $M_B@>\si>>L_B$. If $I\in\ci$ is such that $e_I\in M_B$,
then $\be_I=\p(e_I)\in L_B$ and by 1.2(i) for $\bV$ we have
$I\in B$. We see that $\ph(V)$ satisfies 1.2(i).

For $B\in\ph(V)$ we show:

(a) {\it We have $M_B\sub V_0$ (notation of 3.1). Moreover,
$\p\i(L_B)=M_B\op Fe_S$.}
\nl
By 2.7(a) we can find $t\in S$ such that
$B\in\ph(V(t))$. By 1.5(a) the  subspace of $V$ (or $V(t)$)
spanned by $\{e_I;I\in B\}$ is
contained in $V(t)_0$. Thus, $ M_B\sub V(t)_0$. 

Let $x\in M_B$. We have $x\in V(t)_0$; since
$e_S\n V(t)$ we have  $x=e_I$ for some $I\sub S$, $I\ne S$. By the
definition of $V(t)_0$ we have $|c(I)^{0+}|=|c(I)^{0-}|$ (see 1.5)
so that $|c^0(I)|=|c(I)^{0+}|+|c(I)^{0-}|$ is even and $e_I\in V_0$.
Thus $ x\in V_0$. 
This proves the first assertion of (a).
For the second assertion we note that $M_B$ is a hyperplane
in $\p\i(L_B)$ and that $e_S\in\p\i(L_B)$. It remains to note that
$e_S\n M_B$ (since $e_S\n V(t)$).

\subhead 3.3\endsubhead
Consider the map $\e:\ph(V)@>>>V$ in 1.2(ii). For $B\in\ph(V)$ we show:

(a) {\it We have $\e(B)\in M_B$. In particular we have $\e(B)\in V_0$.}
\nl
(See 3.2(a).)
As in the proof of 3.2(a) we can assume that $B\in\ph(V(t))$
where $t\in S$. Using the commutative diagram in 2.7 we are reduced to
property 1.2(ii) for $V(t)$ which is already known.

\mpb

We show:

(b) {\it The map $\e:\ph(V)@>>>V$ restricts to a bijection
$\ph(V)@>\si>>V_0$.}
\nl
The composition $\p\e:\ph(V)@>>>\bV$ is the same as the map $\e$ for
$\bV$  hence is a bijection. It follows that $\e:\ph(V)@>>>V$ is
injective and its image has exactly $2^{N-1}$ elements. Since this
image is contained in $V_0$ (see (a)) and $|V_0|=2^{N-1}$, we see that
(b) holds.

We show:

(c)  $V_0=\cup_{B\in\ph(V)}M_B$
\nl
The right hand side is contained in the left hand side by 3.2(a).
Now let $x\in V_0$. By \cite{L20} we have 
$\bV=\cup_{L\in\cf(\bV)}L$. Thus, we have $\p(x)\in L_B$ for some
$B\in\ph(V)$. It follows that we have $x\in\p\i(L_B)=M_B\op Fe_S$..
It is enough to show that $x\in M_B$. If $x\n M_B$, then
$x+e_S\in M_B$ so that by (a) we have $x+e_S\in V_0$.
Using 3.1(b) we then have $x\in V_1$, contradicting $x\in V_0$.
This proves (c).

We see that $\ph(V)$ satisfies 1.2(ii),(iii).

Now let $B,B'$ in $\ph(V)$ be such that $\e(B')\in M_B$. Applying
$\p$ we see that $\p\e(B')\in L_B$. Note that $\p\e$ is the same as
$\e$ relative to $\bV$. Since $\ph(\bV)$ satisfies 1.2(iv), we see
that $g_s(B')\le g_s(B)$ for any $s\in S$. (The function $g_s$ is
the same for $V$ as for $\bV$.) Thus, 1.2(iv) holds for $\ph(V)$.
This completes the proof of Theorem 1.4.

\subhead 3.4\endsubhead
Let $B\in\ph(V)=\ph(\bV)$ be such that $B\ne\emp$.
Then $\supp(B)\ne\emp$ and by 2.6 we have $\supp(B)\ne S$ hence the
subset $c(\supp B)$ of $\ci$ is defined as in 1.3. As in the proof of
3.1(b) we have $c(\supp(B))=\{I_1,I_3,\do,I_{2r-1}\}$,
$c(S-\supp(B))=\{I_2,I_4,\do,I_{2r}\}$ for some $r\ge1$.
Since $e_{I_1\cup I_3\cup\do I_{2r-1}}\in V_0$, from 3.1(b) we see that
$e_{I_2\cup I_4\cup\do I_{2r}}\in V_1$, so that

(a) {\it $|I_k|$ is even for some $k\in\{2,4,\do,2r\}$. In particular
there exist $s,s'$ in $S$ such that $\{s,s'\}\in\fE$ and
$\supp(B)\cap\{s,s'\}=\emp$.}
\nl
We show:

(b) $|B|\le(|S|-1)/2$.
\nl
A proof identical to that of \cite{L20, 1.3(g)} shows:

(c) {\it If $I\in B$ then $|\{I'\in B;I'\sub I\}|=(|I|+1)/2$.}
\nl
Using (c) we have
$$\align&|B|=\sum_{I\in c(\supp(B)}
=\sum_{I\in\c(\supp B}|\{I'\in B;I'\sub I\}|\\&
\le\sum_{I\in\c(\supp B}(|I|+1)/2
=(|I_1|+1)/2+(|I_3|+1)/2+\do+(|I_{2r-1}|+1)/2\\&=
(|I_1|+|I_3|+\do+|I_{2r-1}|+r)/2=
(|S|-|I_2|-|I_4|-\do-|I_{2r}|+r)/2\le|S|/2.\endalign$$
Thus $|B|\le|S|/2$. Since $|B|\in\NN$ and $|S|$ is odd we see that (b)
holds.

We show:

(d) {\it We have $|B|=(|S|-1)/2$ if and only if we
have $|I_k|=1$ for all
$k\in\{2,4,\do,2r\}$ except for a single value of $k$ for which
$|I_k|=2$.}
\nl
Assume first that $|B|=(|S|-1)/2$. The proof of (c) shows that in our
case $(|S|-|I_2|-|I_4|-\do-|I_{2r}|+r)/2$ is equal to $(|S|-1)/2$ or
to $|S|/2$, hence $(|I_2|-1)+(|I_4|-1)+\do+(|I_{2r}|-1)$ is equal to
$1$ or $0$. Thus either (d) holds or else
we have $|I_k|=1$ for all $k\in\{2,4,\do,2r\}$ without exception.
This last possibility is excluded by (a). This proves one implication
of (d). The reverse implication follows from the proof of (c).

\subhead 3.5\endsubhead
Let $\ee$ be a two element subset of $S$ such that $\ee\in\fE$.
Let $[\ee]=\be_{(S-\ee)^{odd}}\in\bV$.
We define a linear function $z_\ee:\bV@>>>F$ by
$z_\ee(\be_s)=1$ if $s\in\ee$, $z_\ee(\be_s)=0$ if $s\in S-\ee$.
Note that the radical of $<,>|_{z_\ee\i(0)}$ is $F[\ee]$.

Let $B\in\ph(\bV)$. The following result is used in \cite{L23, 3.5}.

(a) {\it If $[\ee]\in L_B$ then $\supp(B)\cap\ee=\emp$ and
$|B|=(|S|-1)/2$.}
\nl
Let $B^*\in\ph(\bV)$ be the subset of $R$ consisting of the various
$\{s\}$ with $s\in(S-\ee)^{odd}$. We have $[\ee]=\e(B^*)$ so that
$B^*\le B$. Using 1.2(iv), we see that $g_s(B^*)\le g_s(B)$ for all
$s\in S$. It follows that $g_s(B)\ge1$ for all $s\in(S-\ee)^{odd}$.
Thus $(S-\ee)^{odd}\sub\supp(B)$.

Let $\{I_{i_1},I_{i_2},\do,I_{i_l}\}$ be the subset
of $\{I_2,I_4,\do,I_{2r}\}$ consisting of those $I_k$ ($k$ even)
such that $|I_k|\ge2$. This subset is nonempty by 3.4(a).
Let $I\in\{I_{i_1},I_{i_2},\do,I_{i_l}\}$.
We have $I\cap\supp(B)=\emp$ hence $I\cap(S-\ee)^{odd}=\emp$.
If $I\ne\ee$ then, since
$|I|\in\{2,4,6,\do\}$ we have $I\cap(S-\ee)^{odd}\ne\emp$, a
contradiction. Thus, $I=\ee$. We see that
$\ee\cap\supp(B)=\emp$ that is $\supp(B)\sub S-\ee$. Moreover,
$\{I_{i_1},I_{i_2},\do,I_{i_l}\}$ consists of a single object
namely $\ee$. It remains to use 3.4(d).

Conversely,

(b) {\it If $\supp(B)\cap\ee=\emp$ and $|B|=(|S|-1)/2$, then
$[\ee]\in L_B$.}
\nl
Note that $L_B$ is an isotropic
subspace of $\z_\ee\i(0)$ and in fact a maximal one since
$\dim(L_B)=(\dim(\z_\ee\i(0))+1)/2$. But any maximal isotropic subspace
of $\z_\ee\i(0)$ must contain the radical $F[\ee]$. Thus, (b) holds.

\head 4. Complements\endhead
\subhead 4.1\endsubhead
In this subsection we assume that $(V,<>,e:S@>>>V)$ is as in 1.3(a), but
the condition that $N\in\{3,5,7,\do\}$ is replaced by the condition that
$N\in\{4,6,8,\do\}$. From the results in \cite{L19} one can deduce that
$(V,<>,e:S@>>>V)$ is still perfect with $V_0$ having the same
description as in 1.5. 
Let $S'$ be a subset of $S$ such that $S'\in\ci$,
$|S'|=|S|-1$. Let $V'$ be the subspace of $V$ spanned by
$\{e_s;s\in S'\}$. Then $V'$ with the restriction of $<,>$ to $V'$ and
with $S'@>>>V'$, $s\m e_s$ is as in 1.3(a) so that $\ph(V')$ and
the image $V'_0$ of $\e:\ph(V')@>>>V'$ is defined.
Let $S^{odd}\sub S$ be as in 1.1. (This is defined since
$S\in\ci^1$.) Note that the radical of $<,>$ on $V$ is
$Fe_{S^{odd}}$. One can show that

(a) $V_0=V'_0\sqc(V'_0+e_{S^{odd}})$.
\nl
Hence there is a unique fixed point free involution $B\m B'$ of
$\ph(V)$ such that $\e(B')=\e(B)+e_{S^{odd}}$ for all $B\in\ph(V)$.

\subhead 4.2\endsubhead
In this subsection we assume that $(V,<>,e:S@>>>V)$ is as in 1.3(b);
we preserve the notation of \S3.

Let $\cf(V)$ (resp. $\cf^1(V)$)
be the collection of subspaces of $V$ of the form $M_B$ (resp.
$M_B\op Fe_S$) for various $B\in\ph(V)$. Let
$\ti\cf(V)=\cf(V)\sqc\cf^1(V)$.
We show that $\ti\cf(V)$ has properties similar to those  of $\cf(V)$.
We define $\ti\e:\ti\cf(V)@>>>V$ by $\ti\e(M_B)=\e(B)$,
$\ti\e(M_B\op Fe_S)=\e(B)+e_S$.
Note for any $X\in\ti\cf(V)$ we have $\ti\e(X)\in X$.
(This is similar to 1.2(ii).)

Now $\ti\e$ restricts to the bijection $\cf(V)@>\si>>V_0$,
$M_B\m\e(B)$ and to the bijection $\cf^1(V)@>>>V_1$,
$M_B\op Fe_S\m\e(B)+e_S$ (recall the bijection $x\m x+e_S$,
$V_0@>\si>>V_1$). Hence $\ti\e$ is a bijection.
(This is similar to 1.2(iii).)

For $X,X'$ in $\ti\cf(V)$ we say that $X'\le X$ if one of the following
holds:

$X=M_B,X'=M_{B'}$ and $B'\le B$ in the partial order 1.2(a) on $\ph(V)$;

$X=M_B\op Fe_S,X'=M_{B'}\op Fe_S$ and $B'\le B$ in the partial order
1.2(a) on $\ph(V)$;

$X=M_B\op Fe_S,X'=M_{B'}$ and $B'\le B$ in the partial order 1.2(a) on
$\ph(V)$.

This is a partial order on $\ti\cf(V)$. (This is similar to 1.2(iv).)

\subhead 4.3\endsubhead
In this subsection we assume that $(V,<>,e:S@>>>V)$ (as in 1.1)
is perfect. Let $B\in\ph(V)$.
We will give an alternative formula for $\bar\e(B)$.

We define a partition $B=B_1\sqc B_2\sqc B_3\sqc\do$ as follows.

$B_1$ is the set of all $I\in B$ such that $I$ is not properly
contained in any $I'\in B$. Now $B_2$ is the set of all
$I\in B-B_1$ such that $I$ is not properly
contained in any $I'\in B-B_1$. Now $B_3$ is the set of all
$I\in B-(B_1\cup B_2)$ such that $I$ 
is not properly contained in any $I'\in B-(B_1\cup B_2)$, etc.

For $k\ge1$ we set

$v_k(B)=\sum_{I\in B_k}e_I\in V$.
\nl
We have

(a) $\bar\e(B)=v_1(B)+v_3(B)+v_5(B)+\do$.
\nl
Let $s\in S$. There is a unique sequence
$I_1\in B_1,I_2\in B_2,\do,I_l\in B_l$
such that $s\in I_l\sub I_{l-1}\sub\do\sub I_1$ and
$s\notin \cup_{I\in B_{l+1}}I$.
The coefficient of $e_s$ in $v_1(B)+v_3(B)+v_5(B)+\do$
is $0$ if $l=0\mod4$; is $1$ if $l=1\mod4$; is $1$ if $l=2\mod4$;
is $0$ if $l=3\mod4$. We have $g_s(B)=l$. Note that
$(1/2)l(l+1)\mod 2$
is $0$ if $l=0\mod4$; is $1$ if $l=1\mod4$; is $1$ if $l=2\mod4$;
is $0$ if $l=3\mod4$.  This proves (a).

\subhead 4.4\endsubhead
In this subsection we are in the setup of 2.1.
Let $\bV^\CC$ be the $\CC$-vector space of functions
$\bV@>>>\CC$. For any $x\in\bV$ let $f_x\in\bV^\CC$ be the
function which takes value $1$ on the subspace $L_{\bar\e\i(x)}$ of
$\bV$ and the value $0$ on the complement of that subspace;
let $f'_x\in\bV^\CC$ be the
function which takes value $1$ on the subspace
$\{x'\in\bV;<x',L_{\bar\e\i(x)}>=0\}$ of
$\bV$ and the value $0$ on the complement of that subspace.
From Theorem 1.4 we see that for $x\in\bV$ we have
$f'_x=\sum_{y\in\bV}c_{y,x}f_y$ where $c_{y,x}\in\ZZ$.
Moreover, from the triangularity of Fourier transform \cite{L20a}
we see that $c_{y,x}=0$ unless $x=y$ or
$\dim L_{\bar\e\i(x)}<\dim L_{\bar\e\i(y)}$ and that
$c_{x,x}=\pm2^k$ for some $k\in\NN$. We conjecture that

(a) for any $x,y$ in $\bV$, we have either
$c_{y,x}=0$ or $c_{y,x}=\pm 2^k$ for some $k\in\NN$.
\nl
The dihedral group $Di_{2N}$ of order $2N$ acts naturally on $\bV$;
see 1.3. Let $Z_N$ be a set of representatives for the $Di_{2N}$-orbits.
Assume for example that $x=0$. Then $y\m c_{y,0}$ is constant
on each $Di_N$-orbit. We describe this function assuming that
$S=S_N$ (see 2.8) and $N=7$. We can take

(b)  $\{1245\},\{12345\},\{1235\},\{135\},\{123\},\{14\},\{13\},
\{1\},\{\emp\}$
\nl
where we write
$i_1i_2\do i_m$ instead of $\be_{i_1}+\be_{i_2}+\do+\be_{i_m}$.
The value of $y\m c_{y,0}$ at the $9$ elements in (b) (in the order
written) is

$1,0,1,-1,-1,0,1,-2,8$.

\head 5. The set $\o(\bV)$\endhead
\subhead 5.1\endsubhead
In this section we assume that $(\bV,<>,\p e:S@>>>\bV)$ is as in
1.3(c).
We fix a two element subset $\ee$ of $S$ such that $\ee\in\fE$.

\subhead 5.2\endsubhead
For $B\in R$ we set
$$n_B=|\{I\in B;\ee\sub I\}|\in\NN.$$

Let $\ph(\bV)^\ee=\{B\in\ph(\bV);\supp(B)\cap\ee\ne\emp\}$.

If $B\in\ph(\bV)^\ee$ (in particular if $n_B>0$),
then using $(P_0),(P_1)$, we see that
there is a unique $I_B\in B$ such that $|I_B\cap\ee|=1$.

We have $\ph(\bV)^\ee=\sqc_{\t\in\ee}\ph(\bV)^\t$
where $\ph(\bV)^\t=\{B\in\ph(\bV)^\ee;\t\in I_B\}$.

For $B\in\ph(\bV)$ we define $B^!\in R$ by

$B^!=B-\{I_B\}$ if $n_B\in\{1,3,5,\do\}$

$B^!=B$ if $n_B\in\{0,2,4,\}$.

Note that for $B\in\ph(\bV)$ we have $n_{B^!}=n_B$. We show:

(a) {\it If $B\in\ph(\bV)$, $B'\in\ph(\bV)$ satisfy $B^!=B'{}^!$, then
$B=B'$.}
\nl
If $n_B$ is odd, then from the definition we see that $B^!$ does
not satisfy $(P_1)$. Hence to prove (a) we can assume that
both $n_B$ and $n_{B'}$ are odd.

There is a unique $I\in B^!=B'{}^!$ such that $\ee\sub I$
and such that any $I'\in B^!=B'{}^!$ with $I'\prec I$ satisfies
$\ee\cap I=\emp$. We have $I\in B,I\in B'$.
Let $I_1,I_2,\do,I_k$ (resp. $I'_1,I'_2,\do,I'_l$)
be defined in terms of $I$ as in $(P_1)$ for $B$ (resp. $B'$).
We can assume that $I_B=I_1$ (resp. $I_{B'}=I'_1$)
and $I_2,I_3,\do,I_k$
(resp. $I'_2,I'_3,\do,I'_l$) are the maximal objects of $B^!$
(resp. $B'{}^!$) that are strictly contained in $I$.
Hence $\{I_2,I_3,\do,I_k\}=\{I'_2,I'_3,\do,I'_l\}$. Note that
$I_1$ is the unique object of $\ci^1$ such that $I_1\spa I_j$
for $j>1$ and $I^{ev}\sub I_1\sqc I_2\sqc\do\sqc I_k$; similarly
$I'_1$ is the unique object of $\ci^1$ such that $I'_1\spa I'_j$
for $j>1$ (that is $I'_1\spa I_j$ for $j>1$)
and $I^{ev}\sub I'_1\sqc I'_2\sqc\do\sqc I'_l$
(that is $I^{ev}\sub I'_1\sqc I_2\sqc\do\sqc I_k$).
It follows that $I_1=I'_1$ so that $B=B'$. This proves (a).

Let
$$\o(\bV)=\{B^!;B\in\ph(\bV)\}\sub R.$$
From (a) we see that

(b) {\it $B\m B^!$ defines a bijection $\ph(\bV)@>\si>>\o(\bV)$.}
\nl
For any $B\in\o(\bV)$ we define $\tB\in\ph(\bV)$ by $B=\tB^!$.

There is a unique bijection ${}'\e:\o(\bV)@>\si>>\bV$
such that ${}'\e(B)=\bar\e(\tB)$ for any $B\in\o(\bV)$.

\mpb

There is a unique involution $\io:S@>>>S$ preserving the graph
structure and interchanging the two elements of $\ee$.
It induces an involution on $R$ denoted again by $\io$ which leaves
stable $\ph(\bV)$ and $\o(\bV)$.

\subhead 5.3\endsubhead
We now assume that instead of specifying an element $\ee$ of $\fE$
we specify an element $\ee'\in\fE'$ (see 2.8) that is
a pair $\{s_1,s\},\{s_2,s\}$ of two distict two edges of $S$
whose intersection is $\{s\}$ for some $s\in S$.
In terms of $\ee'$ we have a function
$(X_1,X_2,\do,X_k)\m n_{X_1,X_2,\do,X_k}$
from $\ph(\bV')$ (see 2.8) to $\NN$ defined in a way
analogous to the way $B\m n_B$
from $\ph(\bV)$ to $\NN$ was defined in terms of $\ee$. We have 
$$n_{X_1,X_2,\do,X_k}=|\{i\in\{1,2,\do,k\},s\sub\un{X_i}-X_i\}|.$$
The analogue of the assigment $B\m I_B$ for $B\in\ph(\bV)$ such that
$n_B>0$ is the assignent
$$\{X_1,X_2,\do,X_k\}\m I_{\{X_1,X_2,\do,X_k\}}=X$$
for any $\{X_1,X_2,\do,X_k\}\in\ph(\bV')$ such that
$n_{X_1,X_2,\do,X_k}>0$; here $X$ is the unique $X_i$
such that $s\in X_i$.
Then $\o(\bV')$ is defined in terms of $s$ in the same way as $\o(\bV)$
was defined in terms in terms of $\ee$.
Namely $\o(\bV')$ consists of the sequences obtained from various
sequences $\{X_1,X_2,\do,X_k\}\in\ph(\bV')$ by removing
$X=I_{\{X_1,X_2,\do,X_k\}}$ whenever $X$ is defined and by not removing
anything whenever $X$ is not defined.

This approach appears in \cite{L23} (in a less symmetric and more
complicated way) where $S=S_N$ as in 2.8. The set
$\cx_{N-2}$ defined in \cite{L23, 1.3} is the same as $\o(\bV)$ if
$\bV,\bV'$ are identified as in 2.8 and if $\ee$ is taken to be
$\{N-1,N\}$ so that $s=N$.

Hence $\o(\bV)$ is closely related to the theory of
unipotent representations of even orthogonal groups over a finite
field in the same way as $\ph(\bV)$ is closely related to the
theory of unipotent representations of symplectic groups over a
finite field.

\subhead 5.4\endsubhead
For $B\in\o(\bV)$ we denote by $<B>$ the subspace of $\bV$ spanned by
$\{\be_I;I\in B\}$.

For $B',B$ in $\o(\bV)$ we write $B'\preceq B$
if there exists a sequence
$$B'=B_0,B_1,B_2,\do,B_k=B$$
such that
$${}'\e(B_0)\in<B_1>,{}'\e(B_1)\in<B_2>,\do,{}'\e(B_{k-1})\in<B_k>.
\tag a$$
We show:

(b) {\it $\preceq$ is a partial order on $\o(\bV)$.}
\nl
In the setup of (a), for $i=0,1,\do,k$ we have  
$<B_i>\sub L_{\tB_i}$ hence $\bar\e(\tB_i)={}'\e(B_i)\in L_{\tB_i}$.
We see that if $B'\preceq B$ then $\tB'\le \tB$ in
$\ph(\bV)$. It is enough to prove that if $B'\preceq B$ in $\o(\bV)$
and $B\preceq B'$ in $\o(\bV)$ then $B'=B$. We have
$\tB'\le\tB$ in $\ph(\bV)$ and $\tB\le\tB'$ in $\ph(\bV)$.
Since $\le$ is a partial order on $\ph(\bV)$ we have $\tB'=\tB$.
It follows that $B=B'$. This proves (a). (See also
\cite{L23, 2.10(a)}).

\head 6. The subsets $\o^+(bV),\o^-(\bV)$ of $\o(\bV)$\endhead
\subhead 6.1\endsubhead
In this section we preserve the setup of 5.1. Let $z_\ee:\bV@>>>F$
be as in 3.5. Let $\bV^+=z_\ee\i(0),\bV^-=z_\ee\i(1)$.
We set $\o^+(\bV)={}'\e\i(\bV^+),\o^-(\bV)={}'\e\i(\bV^-)$.
We have $\o(\bV)=\o^+(\bV)\sqc\o^-(\bV)$ and ${}'\e$ restricts
to bijections $\o^+(\bV)@>>>\bV^+,\o^-(\bV)@>>>\bV^-$. We show:

(a) {\it If $B\in\ph(\bV), n_B=2k+1$, then $\bar\e(B)\in\bV^+$
so that $B^!\in\o^+(\bV)$.}
\nl
By $(P_1)$ we can find $I'\in B$ such that $I'\cap\ee=\{\s\}$ for
some $\s\in\ee$; let $\s'\in\ee,\s'\ne\s$. We then have $g_\s(B)=2k+2$,
$g_{\s'}(B)=2k+1$. We have
$$\align&\bar\e_\s(B)+\bar\e_{\s'}(B)\\&
=(1/2)(2k+2)(2k+3)+(1/2)(2k+1)(2k+2)\\&=(1/2)(2k+2)(4k+4)=0\mod2
\endalign$$

so that $z_\ee(\bar\e(B))=0$ that is $\bar\e(B)\in\bV^+$.

We show:

(b) {\it If $B\in\ph(\bV), n_B=2k$, $k\ge1$, then $\bar\e(B)\in\bV^-$ so
that $B^!\in\o^-(\bV)$.}
\nl
By $(P_1)$ we can find $I'\in B$
such that $I'\cap\ee=\{\s\}$ for some $\s\in\ee$; let
$\s'\in\ee,\s'\ne\s$. 
We then have $g_\s(B)=2k+1$, $g_{\s'}(B)=2k$. We have
$$\align&\bar\e_\s(B)+\bar\e_{\s'}(B)=\\&=(1/2)(2k+1)(2k+2)
+(1/2)2k(2k+1)\\&=(1/2)(2k+1)(4k+2)=(2k+1)^2=1\mod2\endalign$$

so that $z_\ee(\bar\e(B))=1$ that is $\bar\e(B)\in\bV^-$.
Note that

$\{B\in\o^+(\bV);n_B=0\}=\{B\in\ph(\bV);\supp(B)\cap\ee=\emp\}$,

$\{B\in\o^-(\bV);n_B=0\}=\{B\in\ph(\bV);|\supp(B)\cap\ee|=1\}$.

\subhead 6.2\endsubhead
Let $B'\in\o(\bV)$. We write $B'=B^!$ where $B\in\ph(\bV)$.

Assume first that $B$ is as in 6.1(a). Then $B'\in\o^+(\bV)$ and
$I_B$ is the only $I\in B$ such that
$|I\cap\ee|=1$; since $B^!=B-I_B$ we see that for any $I\in B'$
we have $|I\cap\ee|\in\{0,2\}$.

Assume next that $B$ is as in 6.1(b).
Then $B'=B\in\o^-(\bV)$ and
$I_B$ satisfies $|I_B\cap\ee|=1$; thus, for
some  $I\in B'$ we have $|I\cap\ee|=1$,

We now assume that $n_B=0$. If $\supp(B)\cap\ee=\emp$, then clearly 
we have $|I\cap\ee|=0$ for any $I\in B$.
If $|\supp(B)\cap\ee|=1$, then clearly
we have $|I\cap\ee|=1$ for some $I\in B$.

We see that for $B\in\o(\bV)$ the following holds:

(a) {\it $B\in\o^+(\bV)$ if and only if
$|I\cap\ee|\in\{0,2\}$ for any $I\in B$.}

\subhead 6.3\endsubhead
We show:

(a) {\it Let $B',B$ in $\o(\bV)$ be such that $B'\preceq B$. If
$B\in\o^+(\bV)$, then $B'\in\o^+(\bV)$.}
\nl
We can assume that ${}'\e(B')\sub<B>$. (The general case would
follow by using several times this special case.) By 6.2(a) we
have $|I\cap\ee|\in\{0,2\}$ for any $I\in B$. It follows that any
$x\in<B>$ satisfies $z_\ee(x)=0$. In particular we have
$z_\ee({}'\e(B'))=0$ so that ${}'\e(B')\in\bV^+)=0$ and
$B'\in \o^+(\bV)$. This proves (a).

\define\bbV{\ov{\bV}}

\head 7. The sets $\cf^+(\bbV)^\t,\cf^-(\bbV)^\t$ \endhead
\subhead 7.1\endsubhead
In this section we preserve the setup of 5.1.  
For $\t\in\ee$ let $\o(\bV)^\t=\{B\in\o(\bV);\tB\in\ph(\bV)^\t\}$.
We have
$\o(\bV)^\t=\o^+(\bV)^\t\sqc\o^-(\bV)^\t$
where for $\d\in\{+.-\}$ we set$\o^\d(\bV)\t=\o(\bV)^\t\cap\o^\d(\bV)$.

Under the identification $\o(\bV)=\o(\bV')$ in 2.8, 5.3 and with
notation of \cite{L23, 1.4}, the following holds:

If $n\in\{1,3,5,\do\}$, then

$\{B\in\o^+(\bV)^{N-1},n_B=n\}$ becomes $\cx_{N-2}^{t,+}$, $t=-n-1$;

$\{B\in\o^+(bV)^{N},n_B=n\}$ becomes $\cx_{N-2}^{t,+}$, $t=n+1$;   

if $n\in\{0,2,4,6,\do\}$, then

$\{B\in\o^-(\bV)^{N-1},n_B=n\}$ becomes $\cx_{N-2}^{t,-}$, $t=n$;

$\{B\in\o^-(\bV)^{N},n_B=n\}$ becomes $\cx_{N-2}^{t,-}$, $t=-n-2$.

\subhead 7.2\endsubhead
Let $\t\in\ee$.

(a) {\it Assume that $B'\in\o^+(\bV),B\in\o^+(\bV)^\t$ satisfy
$B'\preceq B$ and $n_B>0$. Then we have either $n_{B'}=n_B$ and
$B'\in\o^+(\bV)^\t$, or else $n_{B'}<n_B$.}

(b) {\it Assume that $B'\in\o^-(\bV),B\in\o^-(\bV)^\t$ satisfy
$B'\preceq B$ and $n_B\ge0$. Then we have either $n_{B'}=n_B$ and
$B'\in\o^-(\bV)^\t$, or else $n_{B'}<n_B$.}
\nl
Using the identification $\o(\bV)=\o(\bV')$ in 2.8, 5.3 and the
results in 7.1 we see that when $\t=N-1$, (a) follows from
\cite{L23, 3.2} and (b) follows from \cite{L23, 3.4}. Using the
symmetry $\io$, we see that (a) and (b) for $\t=N$ follow from (a)
and (b) for $\t=N-1$.

\subhead 7.3\endsubhead
We choose a subset $J$ of $S-\ee$ such that $|J|=N-3$ and
such that when $N>3$ we have $J\sub\ci$.

Let $\o(\bV)_J=\{B\in\o(\bV);\supp B\sub J\}$. Then ${}'\e$ defines a
bijection of $\o(\bV)_J$ onto a subset $\bV_{J,0}$ of $\bV$. We set
$$\bV_{J,1}={}'\e(\{B\in\o(\bV);\supp(B)\cap\ee=\emp\})-\bV_{J,0}
\sub\bV.$$

Assume now that $B'\in\o(\bV),B\in\o(\bV)_J$
satisfy $B'\preceq B$. From \cite{L23, 3.3} we deduce:

(a) {\it We have $B'\in\o(\bV)_J$.}

\subhead 7.4\endsubhead
Let $\t\in\ee$. We set $\ti\o^+(\bV)^\t=\o^+(\bV)^\t\cup\o(\bV_J)$
$\ti\o^-(\bV)^\t=\o^-(\bV)^\t$. 

Assume now that $B'\in\o^\d(\bV),B\in\ti\o^\d(\bV)^\t$ satisfy
$B'\preceq B$. From 7.2(a),(b) and 7.3(a) we deduce:

(a) {\it We have either $B'\in\ti\o^\d(\bV)^\t$ and $n_{B'}=n_B$, or else
$n_{B'}<n_B$.}

\subhead 7.5 \endsubhead
Let $\bbV=\bV/F[\ee]$ and let $\bp:\bV@>>>\bbV$ be the obvious
quotient map.
Let $\bbV^+=\bp(\bV^+)$, $\bbV^-=\bp(\bV^-)$.
We have $[\ee]\in\bV^+$ hence $\bbV=\bbV^+\sqc\bbV^-$
and $|\bbV^+|=(1/2)|\bV^+|=|\bbV^-|$.

Let $\d\in\{+,\}$. For $n\ge0$, $\t\in\ee$ we set
$$\bV^{\d,\t}_n={}'\e(\{B\in\o^\d(\bV)^\t;n_B=n\})\sub\bV^\d.$$
From the results in \cite{L23, 2.7, 3.5} we see that

(a) {\it the two subsets $\bV^{\d,\t}_n$ (with $\t\in\ee$) are interchanged by 
the involution $x\m x+[\ee]$ of $\bV^\d$;}

(b) {\it $\bV_{J,0},\bV_{J,1}$ are interchanged by the involution
$x\m x+[\ee]$ of $\bV$.}
\nl
(For (b) see also 4.1(a).)

For $\t\in\ee$ we set

$$H^{\d,\t}={}'\e(\ti\o^\d(\bV)^\t)\sub\bV^\d$$.
We have
$$H^{+,\t}=\bV_{J,0}\cup\cup_{n\ge0}\bV^{+,\t}_n,$$
$$H^{-,\t}=\cup_{n\ge0}\bV^{-,\t}_n$$
From (a),(b) we see that $\bp$
restricts to  bijections $H^{\d,\t}@>\si>>\bbV^\d$.

For $y\in\bbV^\d$ we denote by $\ty^\t\in H^{\d,\t}$ the inverse image
of $y$ under this bijection and we define $\nu_y\in\NN$ by:

$\nu_y=n$ if $\ty^\t\in\bV^{\d,\t}_n$,

$\nu_y=0$ if $\d=+$ and $\ty^\t\in\bV_{J,0}$.

\subhead 7.6 \endsubhead
Let $\d\in\{+.-\}, \t\in\ee$.
For $y',y$ in $\bbV^\d$ we say that $y'\le_\t y$ if there exists

(a) {\it a sequence $y'=y_0,y_1,y_2,\do,y_k=y$ in $\bbV^\d$ such that for
$i\in\{0,1,\do,k-1\}$ we have $\ty_i^\t\in<{}'\e\i(\ty_{i+1}^\t)>$
or $\ty_i^\d+[\ee]\in<{}'\e\i(\ty_{i+1}^\d)>$.}
\nl
We show that in this situation, for any $i\in\{0,1,\do,k-1\}$ we have

(b) $\nu_{y_i}\le\nu_{y_{i+1}}$.
\nl
We set $B_i={}'\e\i(\ty_i^\t)$, $B'_i={}'\e\i(\ty_i^\t+[\ee])$,
$B_{i+1}={}'\e\i(\ty_{i+1}^\t)$.

If $\ty_i^\t\in<{}'\e\i(\ty_{i+1}^\t)>$, then $B_i\preceq B_{i+1}$ so that by
7.4(a) we have $n_{B_i}\le n_{B_{i+1}}$. But $n_{B_i}=\nu_{y_i}$,
$n_{B_{i+1}}=\nu_{y_{i+1}}$, so that (b) holds.

If $\ty_i^\t+[\ee]\in<{}'\e\i(\ty_{i+1}^\t)>$, then $B'_i\preceq B_{i+1}$,
so that by 7.4(a) we have $n_{B'_i}\le n_{B_{i+1}}$. But
$n_{B'_i}=\nu_{y_i}$, $n_{B_{i+1}}=\nu_{y_{i+1}}$, so that (b) holds.

We now see:

(c) {\it If $y'\le_\t y$, then $\nu_{y'}\le\nu_y$.}
\nl
We show:

(d) {\it $\le_\t$ is a partial order on $\bbV^\d$.}
\nl
For $y\in\bbV^\d$ we have $\ty^\t\in<{}'\e\i(\ty^\t)>$ so that $y\le_\t y$. It
remains to show that

(e) {\it if $y,y'$ in $\bbV^\d$ satisfy $y\le_\t y'$ and $y'\le_\t y$, then $y=y'$.}
\nl
Using (c) we have $\nu_{y'}\le\nu_y$ and $\nu_y\le\nu_{y'}$, hence
$\nu_y=\nu_{y'}$. Consider now a sequence $y'=y_0,y_1,y_2,\do,y_k=y$
as in (a). Using (b) and $\nu_y=\nu_{y'}$ we see that for
$i\in\{0,1,\do,k-1\}$ we have  $\nu_{y_i}=\nu_{y_{i+1}}$.
Recall that we have either

(i)  $B_i\preceq B_{i+1}$, or

(ii)  $B'_i\preceq B_{i+1}$,
\nl
where as before we set $B_i={}'\e\i(\ty_i^\t)$, $B'_i={}'\e\i(\ty_i^\t+[\ee])$,
$B_{i+1}={}'\e\i(\ty_{i+1}^\t)$. Note that $n_{B_i}=n_{B'_i}=n_{B_{i+1}}$.

We have $B_i\in\ti\o^\d(\bV)^\t$, $B'_i\in\ti\o^\d(\bV)^{\t'}$,
$B_{i+1}\in \ti\o^\d(\bV)^\t$, where $\t'\in\ee$ and $\t\ne\t'$.
Using 7.4(a), we see that if (ii) holds, then (since $n_{B'_i}=n_{B_{i+1}}$) we would have $\t=\t'$, a
contradiction. Thus, (i) holds. Using this for $i=0,1,\do,k-1$ we see
that
$$B_0\preceq B_1\preceq B_2\le\do\preceq B_k.$$
In particular we have $B'\preceq B$.
Reversing the roles of $y,y'$ we have similarly $B\preceq B'$.
Since $\preceq$ is a partial order on $\o(\bV)$, it follows that
$B=B'$. Applying ${}'\e$, we obtain $\ty^\t=\ty'{}^t$ hence $y=y'$. This
proves (e) and hence (d).

\subhead 7.7\endsubhead
Let $\d\in\{+.-\}, \t\in\ee$.
For any $y\in\bbV^\d$ we set $<y>_\t:=\bp(<{}'\e\i(\ty^\t)>)$ (a subspace
of $\bbV$) and $<y>_{\t,\d}=<y>_\t\cap\bbV^\d$. Note that if $\d=+$ then
$<y>_{\t,\d}=<y>_\t$; if $\d=-$ then $<y>_{\t,\d}$ is the complement in $<y>_\t$
of a hyperplane of $<y>_\t$. Now, the condition that

$\ty_i^\d\in<{}'\e\i(\ty_{i+1}^\t)>$
or $\ty_i^\d+[\ee]\in<{}'\e\i(\ty_{i+1}^\d)>$
\nl
(in 7.6(a))
is equivalent to the condition that $y_i\in\bp(<{}'\e\i(\ty_{i+1}^\t)>)$.
Thus, the condition that $y,y'$ in $\bbV^\d$ satisfy $y'\le_\t y$ is
equivalent to the following condition:

there exists a sequence $y'=y_0,y_1,y_2,\do,y_k=y$ in $\bbV^\d$ such
that for $i\in\{0,1,\do,k-1\}$ we have $y_i\in<y_{i+1}>_{\t,\d}$.

Let $\cf^\d(\bbV)^\t$ be the collection of subsets of $\bbV^\d$ of the
form $<y>_{\t,\d}$ for various $y\in\bbV^\d$. We show:

(a) {\it If $y',y$ in $\bbV^\d$ satisfy $<y'>_{\t,\d}=<y>_{\t,\d}$, then $y=y'$.}
\nl
Indeed, we have $y\in<y>_{\t,\d}$, $y'\in<y'>_{\t,\d}$, hence
 $y\in<y'>_{\t,\d}$, $y'\in<y>_{\t,\d}$, so that $y\le_\t y'$, $y'\le_\t y$.
Since $\le_\t$ is a partial order, it follows that $y=y'$, proving (a).

We show:

(b) {\it The map $\ti\o^\d(\bV)^\t@>>>\cf^\d(\bbV)^\t$, ${}'\e\i(\ty^\t)\m<y>_{\t,\d}$
(for $y\in\bbV^\d$) is bijective.}
\nl
This map is obviously surjective. Moreover we have
$|\ti\o^\d(\bV)^\t@>>>\cf^\d(\bbV)^\t|=|\bbV^\d|$.
It is then enough to show that $|\cf^\d(\bbV)^\t|=|\bbV^\d|$. This follows
from (a).

We show:

(c) {\it If $y\in\bbV^\d$ and $B={}'\e\i(\ty^\t)$ so that $<y>_{\t,\d}=\p(<B>)$
then $\bp$ restricts to an isomorphism $<B>@>\si>><y>_{\t,\d}$.}
\nl
Indeed it is enough to show that $[\ee]\n<B>$. But in fact we have
even $[\ee]\n L_B$ as a consequence of 3.5(a).

\subhead 7.8\endsubhead
Now the two sets $\cf^-(\bbV)^\t$ (for the two values of $\t\in\ee$)
are interchanged by the involution induced by $\io$; they do not
depend on the choice of $J$ in 7.3.
This is not so for the two sets $\cf^+(\bbV)^\t$ (for the two values of $\t\in\ee$), at least if $N>3$; these sets do depend on the
choice of $J$ in 7.3.
But we prefer one of them over the other; namely we prefer the
value of $\t$ such that $\t$ is not joined in our graph to any
element of $J$. (This determines $\t$ uniquely if $N>3$.)
This is the choice made in \cite{L23}.

\enddocument